\documentstyle[floats,aps]{revtex}
\advance\topmargin 30pt
\begin{document}
\input psfig
\pssilent
\title{A spin network generalization of the Jones Polynomial and
Vassiliev invariants}

\author{Rodolfo Gambini$^{1}$\footnote{Associate member of ICTP.}, 
Jorge Griego$^1$, Jorge Pullin$^2$}
\address{1. Instituto de F\'{\i}sica, Facultad de Ciencias, 
Tristan Narvaja 1674, Montevideo, Uruguay}
\address{
2. Center for Gravitational Physics and Geometry, Department of
Physics,\\
The Pennsylvania State University, 
104 Davey Lab, University Park, PA 16802}
\date{Nov 8th, 1997}
\maketitle
\begin{abstract}
We apply  the ideas of Alvarez and Labastida to the invariant of 
spin networks defined by Witten and Martin based on Chern--Simons theory.
We show that it is possible to define ambient invariants of spin networks
that (for the case of $SU(2)$) can be considered as extensions to 
spin networks of the Jones polynomial. Expansions of the coefficients of 
the polynomial yield primitive Vassiliev invariants. The resulting 
invariants are candidates for solutions of the Wheeler--DeWitt equations
in the spin network representation of quantum gravity.
\end{abstract}

\vspace{-8.5cm} 
\begin{flushright}
\baselineskip=15pt
CGPG-97/10-4  \\
q-alg/9711014\\
\end{flushright}
\vspace{8cm}

Witten \cite{Wi88} showed that the expectation value of a Wilson loop
in an $SU(2)$ Chern--Simons theory,
\begin{equation}
<W(\gamma)> = \int DA\, W(\gamma,A) \exp\left({{i k\over 4\pi} 
\int d^3x {\rm Tr}(A\wedge \partial A +{2\over 3}
A \wedge A \wedge A)}\right)
\end{equation}
is a knot invariant closely related with the Jones polynomial.  More
precisely, it is given by the Kauffman bracket
\cite{Kaknph} knot polynomial. The latter differs from the Jones
polynomial by an overall ``phase factor'' given by the exponential of
the writhe of the knot. This factor accounts for the fact that the
Kauffman bracket is a regular isotopy invariant (it is not invariant
when one removes twists, i.e., it is an invariant of framed knots)
whereas the Jones polynomial is an ambient isotopic invariant.
The explicit relation between the Kauffman bracket and the expectation
value we introduced above is given by $K(\gamma,q)= <W(\gamma)>$ with
$q=\exp(2\pi i /k)$, where $q$ is the variable of the polynomial and
$k$ is the coupling constant of the Chern--Simons
theory\footnote{There has been some controversy \cite{Gu} concerning
the relationship of $q$ and $k$. Nonperturbative effects seem to 
fix the value to be $q=\exp(2\pi i /{k+2})$ \cite{Wi88,Aw}. The
constructions we will present go through for either relationship.}.

The expression of the invariants one obtains from the path integral
appear as power series expansions in terms of the parameter ${2 \pi
i/k}$. It is tantamount to evaluating the polynomials for $q=\exp(2\pi
i /k)$ and expanding the exponential in power series. If one does this
for the Jones polynomial, the resulting coefficients are Vassiliev
invariants \cite{BaNa,BiLi,Ba}.

In a separate context, it was also noticed that the resulting
Vassiliev invariants are candidates for wavefunctions in the loop
representation of quantum gravity
\cite{BrGaPuprl,BrGaPunpb,BrGaPugrg,GaPubook}. In this context,
however, one is not interested in arbitrary knot invariants, but
actually in those that arise as counterparts of states in the
``connection representation'',
\begin{equation}
\psi(\gamma) = \int DA W(\gamma,A) \psi(A).
\end{equation}
This poses a series of restrictions on the possible knot invariants
$\psi(\gamma)$ that one could consider as wavefunctions of a quantum
theory in the loop representation. These restrictions are known as
``Mandelstam identities''. If one concentrates of wavefunctions 
depending on a single loop, they are linear, loop dependent, relations
$\sum_i c_i \psi(\gamma_i)=0$. They arise from the fact that the
Wilson loop $W(\gamma,A)$ is a trace of a group element \cite{GaTr}. 
Of all the Vassiliev invariants that one constructs using the 
techniques described above, we only know for a fact that two of them
are compatible with these relations. This makes it particularly
difficult to construct candidates for quantum states of the
gravitational field in terms of these invariants. Implicit in this
discussion is the fact that in quantum gravity one is interested in 
loops with generic intersections, otherwise the spectra of various
quantum operators is not appropriately realized \cite{BrPu91}. Most of
the discussions of knot invariants stemming from Chern--Simons
theories have concentrated on smooth non-intersecting loops or
multiloops. 

Given the limitations imposed by the Mandelstam identities and their
complicated nature, it has proven very fruitful to consider holonomies
in all possible representations of the group. In a sense this is
tantamount to dealing with ``all possible loops, multiloops and
intersections'' in a unified framework \cite{RoSm}. This leads to the
construction of generalized gauge invariant quantities dependent on
trivalent (or higher) graphs \cite{Wi89,RoSm,GaGrPu}. These invariants
can be thought of as ``Wilson nets'' (generalizations of the Wilson
loops dependent on a spin network), and one can therefore consider
computing their expectation value in a Chern--Simons theory. This was
first attempted by Witten \cite{Wi89} and further developed by Martin
\cite{Ma}. These invariants are related to families of invariants 
constructed using braid group techniques by Kauffman and Lins
\cite{KaLi}. In a recent paper we have considered related invariants as
possible candidates for states of the quantum gravitational field
\cite{GaGrPu}. Unfortunately, a significant limitation of these
invariants is that they are regular isotopic invariants, that is, they
are framing-dependent. Therefore they cannot be genuinely considered
quantum states of the gravitational field, since the latter are
supposed to be diffeomorphism invariant, that is, unchanged by any
smooth deformation of the loops (including twists).

The purpose of this note is to address the latter issue. We will make
use of a series of observations by Alvarez and Labastida \cite{AlLa}
that will allow us to construct an ambient isotopic invariant from the
expectation value of a Wilson-net. The resulting invariant will be
related to the Jones polynomial and will allow us to introduce
Vassiliev invariants in the context of spin nets. These
invariants are therefore candidates for wavefunctions of the quantum
gravitational field. 

Returning to the single-loop (ie, not spin-net) context, it was noted
by Alvarez and Labastida\cite{AlLa}, that one could re-express the
expectation value of the Wilson loop as,
\begin{equation}
<W(\gamma)> = K_0(\gamma) \exp\left(\sum_{i=1}^\infty \sum_{j=1}^{d_i}
\alpha_{ij}(\gamma) r_{ij} \left[{2\pi i \over k}\right]^i\right),
\end{equation}
where $r_{ij}$ are group-dependent constants (the construction goes
through for any group) and $\alpha_{ij}(\gamma)$ are the ``primitive''
Vassiliev invariants\footnote{With the addition of $\alpha_{11}$, which is the
framing-dependent factor, which would be the only contribution to the
expectation value if one considered an Abelian theory.}. These are the
minimal set of invariants for a given order that are needed to
generate all invariants of higher orders via products. What has been
accomplished here is a re-summation of the power series that arises in
computing the expectation value, into an exponential involving a
linear combination of the primitive Vassiliev invariants. This
resummation can be straightforwardly understood in terms of Feynman
diagrammatics, through the usual procedure to generate connected
diagrams by writing the generating function as an exponential
\cite{It}.

Based on this resummation, it is immediate to define an ambient
isotopic invariant for multiloops $P(\gamma,k)$, simply by removing
the framing-dependent factor explicitly, i.e.,
\begin{eqnarray}
<W(\gamma)> &=& \exp\left(\alpha_{11} r_{11} \left[{2\pi i \over
k}\right] \right)   K_0(\gamma)  P(\gamma,k)\label{labastida}\\
P(\gamma,k) &=& \exp\left(\sum_{i=2}^\infty 
v_{i}(\gamma) \left[{2\pi i \over k}\right]^i\right),\\
v_{i} &=& 
\sum_{j=1}^{d_i} \alpha_{ij}(\gamma) r_{ij}.
\end{eqnarray}
This expression very clearly embodies the different properties of
transformation of the invariant. The exponential prefactor is
responsible for the framing dependence (non-invariance under
twists). The $K_0(\gamma)$ factor takes care of the fact that the
invariant changes by a power of $(-1)$ when one ``flips over'' an
entire knot through a diffeomorphism, since the orientation of all
trivalent intersections is reversed (see below).

What we will do in this paper is to follow a similar procedure in the
context of spin networks for the $SU(2)$ group. What allows us to do
this is that the formulae (\ref{labastida}) have been shown to hold 
for multi-loops \cite{AlLa}. We will assume they hold in the case of
intersecting loops, and therefore for spin-nets. This will allow us to
define an ambient isotopic invariant associated with Chern--Simons
theory. In order to do this, we start by recalling the invariant
constructed by Witten and Martin \cite{Wi89,Ma}. It is defined by the
following skein relations,
\begin{eqnarray}
\Delta_j &=&E\left(\raisebox{-5mm}{\psfig{figure=unknotj.eps,height=10mm}}
\,,\,k\right) = {q^{j+{1\over 2}} - q^{-j-{1\over 2}}\over q^{1\over
2} - q^{-{1\over 2}}} \label{delta}\\
\Theta(j_1,j_2,j) &=& E\left(\raisebox{-12.5mm}
{\psfig{figure=thetajj1j2.eps,height=25mm}}\,,\,k\right) = 
\sqrt{\Delta_1 \Delta_2 \Delta_3}\label{theta}\\
E\left(\raisebox{-10mm}{\psfig{figure=y.eps,height=20mm}}\,,\,k\right)
&=& (-1)^{j_1+j_2+j_3}\exp\left(i\pi (h_1+h_2-h_3)\right) 
E\left(\raisebox{-10mm}{\psfig{figure=ytwist.eps,height=20mm}}\,,\,k\right),
\qquad h_i \equiv {j_i (j_i+1)\over k}\label{ytwist}\\
E\left(\raisebox{-10mm}{\psfig{figure=j.eps,height=20mm}}\hspace{-1cm},k\right)
&=& \exp\left(- 2 \pi i h_j\right)
E\left(\raisebox{-10mm}{\psfig{figure=jrtwist.eps,height=20mm}}
\hspace{-0.4cm},k\right)
=
\exp\left( 2 \pi i h_j\right)
E\left(\raisebox{-10mm}{\psfig{figure=jltwist.eps,height=20mm}}
\hspace{-0.4cm},k\right)\label{twist}\\
E\left(\raisebox{-10mm}{\psfig{figure=bola.eps,height=20mm}}
\hspace{-0.8cm},k\right)
&=& {\delta_{j\,j'}\over \Delta_j}
E\left(\raisebox{-10mm}{\psfig{figure=bolaclosed.eps,height=20mm}}\,,\,k\right)
E\left(\raisebox{-10mm}{\psfig{figure=j.eps,height=20mm}}
\hspace{-1cm},k\right)\\
E\left(\raisebox{-10mm}{\psfig{figure=bola3.eps,height=20mm}}
,k\right)
&=& {1 \over \sqrt{\Delta_1 \Delta_2 \Delta_3}}
E\left(\raisebox{-10mm}{\psfig{figure=bola3closed.eps,height=20mm}}\,,
\,k\right)
E\left(\raisebox{-10mm}{\psfig{figure=tri.eps,height=20mm}},k\right)\\
E\left(\raisebox{-10mm}{\psfig{figure=para.eps,height=20mm}}
,k\right)
&=& \sum_{i=|j_1-j_2|}^{i=|j_1+j_2|}\sqrt{\Delta_i \over \Delta_1 \Delta_2}
E\left(\raisebox{-10mm}{\psfig{figure=parawish.eps,height=20mm}}\,,
\,k\right)\\
E\left(\raisebox{-10mm}{\psfig{figure=cruc.eps,height=20mm}}
,k\right)
&=& \sum_{i=|j_1-j_2|}^{i=|j_1+j_2|}
\sqrt{\Delta_i\over \Delta_1 \Delta_2} (-1)^{j_1+j_2+j_3} 
\exp\left(i \pi (h_1+h_2+h_3)\right)
E\left(\raisebox{-10mm}{\psfig{figure=parawish.eps,height=20mm}}\,,
\,k\right)\\
E\left(\raisebox{-10mm}{\psfig{figure=wishbone.eps,height=20mm}}\,,\,
k\right)&=&\sum_{j={|j_1-j_4|}}^{|j_1+j_4|}
\left\{
\begin{array}{ccc}
j_2&j_1&j\\
j_4&j_3&l
\end{array}
\right\}_q 
E\left(\raisebox{-10mm}{\psfig{figure=doubley.eps,height=20mm}}\,,\,k\right),
\label{recoupling}
\end{eqnarray}
where the expression in curly braces is the q-deformed Racah symbol,
which is defined as,
\begin{equation}
\left\{
\begin{array}{ccc}
j_2&j_1&j\\
j_4&j_3&l
\end{array}
\right\}_q 
= {E\left(
\raisebox{-10mm}{\psfig{figure=tet.eps,height=20mm}}\,,\,k\right) 
\over \sqrt{\Delta_1 \Delta_2 \Delta_3 \Delta_4}}.
\end{equation}

The above expressions completely characterize the invariant for any 
spin network. The explicit expression for the value of the invariant
in the last expression ``the tetrahedron diagram'' can be computed
using the recoupling identity (\ref{recoupling}) in conjunction with
the definitions of the Delta and Theta diagrams
(\ref{delta},\ref{theta}) (see \cite{Ma} for an explicit computation).

This invariant is closely related to that of Kauffman and Lins
\cite{KaLi}, but has two main differences. The first one is the choice
in the definition of the ``Wilson net'' gauge invariant object. These
objects are defined by joining holonomies along the lines of the 
net at the vertices using group invariant tensors. There is an 
ambiguity in the definition, given by which line ends up
being contracted into which entry of the $3jm$ coefficients (closely
related to the Clebsch-Gordan coefficients) \cite{VaMoKh}.
This ambiguity can be coded into an ``orientation'' of the lines
entering the intersection, due to the cyclic property of the
coefficients. If one performs a diffeomorphism on a diagram, one 
in general can reverse the orientation of an intersection. Witten and
Martin \cite{Wi89,Ma} introduce an additional normalization factor in 
the definition of the Wilson net in such a way as to make them
invariant under changes of orientation at the vertices. The
coefficient is,
\begin{equation}
V_\pm = \exp\left(\pm {i \pi\over 2} [j_1+j_2+j_3]\right) \sqrt[4]{
(2 j_1+1)(2 j_2+1) (2 j_3+1)}
\end{equation}
where the $\pm$ refers to the two possible orientations of a trivalent
vertex. Kauffman and Lins do not introduce this normalization factor.
The invariant they construct is therefore less naturally associated
with an expectation value of a Wilson net (their construction is
entirely in terms of braid group techniques, and does not concern
itself with viewing the invariant as an expectation value of a Wilson
net). From our point of view, since we are interested in constructing
states of quantum gravity using these invariants, it is therefore
natural to choose the normalization of Witten and Martin.

The second difference between the Witten--Martin and the
Kauffman--Lins invariants is given by the definition of the planar
diagrams (\ref{delta},\ref{theta}). This difference can be understood
\cite{Wi89} in terms of the freedom that one has in choosing the
measure to define the path integral, and the two definitions can be
viewed as a choice of measure. It is not clear at the moment which is
the ``correct'' measure for quantum gravity. Although there is
progress in understanding diffeomorphism invariant measures
\cite{AsLe}, these have proved ineffective at the time of performing
integrals like the expectation value of the Wilson net in a
Chern--Simons theory. Variational techniques \cite{GaPucmp,GaGrPu},
which only involve weak assumptions about the measure considered for
evaluating the path integral, yield unique results for the identities
involving crossings and twists, but up to now they have failed to
produce unique results for the recoupling identities and the
evaluation of the planar diagrams.

We now proceed to apply the technique of Alvarez and Labastida to the
invariant we just introduced. The first step is to understand the
moves that are associated with diffeomorphisms, since one would expect
the final object to be invariant under diffeomorphisms. The moves in 
equations (\ref{ytwist}) and (\ref{twist}) are associated with 
diffeomorphisms. In addition to these, one has to be concerned with
twists at intersections, similar to those in equation (\ref{ytwist})
but involving several strands, specifically, 
\begin{eqnarray}
E\left(\raisebox{-10mm}{\psfig{figure=h.eps,height=20mm}}\,,\,k\right)
&\rightarrow& 
E\left(\raisebox{-10mm}{\psfig{figure=htwist.eps,height=20mm}}\,,\,k\right),
\label{twist2}\\
E\left(\raisebox{-10mm}{\psfig{figure=h2.eps,height=20mm}}\,,\,k\right)
&\rightarrow& 
E\left(\raisebox{-10mm}{\psfig{figure=h2twist.eps,height=20mm}}\,,\,k\right),
\label{twist3}
\end{eqnarray}
and similarly for $n$ insertions.

The first task consists in identifying the equivalent of $v_1\equiv
r_{11} \alpha_{11}$
in  equation (\ref{labastida}), which we can easily do by considering 
the expansion of $E(\Gamma,k)$ to first order in 
$\kappa ={2\pi i\over k}$ 
\begin{equation}
E(\Gamma,\kappa) = E(\Gamma,0) + \kappa E(\Gamma,0) v_1(\Gamma)+\cdots
\end{equation} 
If one now considers equation (\ref{twist}) and its expansion to first
order, it is immediate to see that,
\begin{equation}
v_1\left(\raisebox{-10mm}{\psfig{figure=j.eps,height=20mm}}\hspace{-1cm}\right)
= -j(j+1)+
v_1\left(\raisebox{-10mm}{\psfig{figure=jrtwist.eps,height=20mm}}
\hspace{-0.4cm}\right)
=j(j+1)+
v_1\left(\raisebox{-10mm}{\psfig{figure=jltwist.eps,height=20mm}}
\hspace{-0.4cm}\right).
\end{equation}

What we now see is that if one defines an invariant $P$ as we
introduced in equation (\ref{labastida}), the skein relations we just
discussed for $v_1$ compensate exactly the changes that the invariant
$E$ underwent under twists as in equation (\ref{twist}). It is
straightforward to see that this is also true for twists at the
intersections as those shown in equation (\ref{ytwist}). The only task
left is to show that this is also true for twists with many insertions
of lines, as those shown in (\ref{twist2},\ref{twist3}) and their
generalization for $n$ lines. The way to show that these twists can be
``undone'' via the application of the fundamental Reidemeister moves
and the invariance under the twist at a trivalent vertex we have just
discussed. Specifically, let us consider the following move,
\begin{equation}
E\left(\raisebox{-10mm}{\psfig{figure=htwist.eps,height=20mm}}\,,\,k\right)
\rightarrow 
E\left(\raisebox{-10mm}{\psfig{figure=htwistp.eps,height=20mm}}\,,\,k\right)
=
E\left(\raisebox{-10mm}{\psfig{figure=h.eps,height=20mm}}\,,\,k\right).
\end{equation}

What occurred here is that line number 1 was ``slid'' under the rest
of the lines towards the bottom. This is a valid Reidemeister move as
long as the orientation of the trivalent intersection of the left is
kept unchanged. Once the line number 1 slides below the horizontal
line, one is introducing a twist of the form (\ref{twist}). We have
just shown that the invariant does not change under such twists. From
there on one keeps on performing a Reidemeister move. A similar
discussion applies to the line number 2. We end up with a diagram 
topologically equivalent to the untwisted one. This discussion can be
easily generalized to the case of twists involving $n$ inserted
lines. 

We have therefore constructed an invariant polynomial which admits an
expansion in terms of the analog for spin-networks of Vassiliev
invariants. Functions of these invariants are natural candidates for
states of the quantum gravitational field in the loop representation
of quantum gravity. In particular, they are diffeomorphism invariant,
compatible with the Mandelstam identities of $SU(2)$ and particularize
for the case of single loops to invariants that in some particular
cases have been found to be annihilated by the Hamiltonian constraint
of quantum gravity. 

It has been conjectured \cite{BaNa} that Vassiliev invariants could be
able to separate all knots. The basis of this conjecture is that it is
true for braids. This suggests that any knot invariant could be
expanded in terms of these invariants. If one believes this
conjecture, one should conclude that these invariants would constitute
an ideal setting for discussing quantum gravity in the loop
representation. An important issue to be addressed, as pointed out by
Kauffman \cite{KaBa} is if all Vassiliev invariants arise from
Chern--Simons theory. If this were the case, then the actual setting
for quantum gravity could be given completely by the invariants we
have been discussing.  It is now known that all Vassiliev invariants
can be obtained from Chern--Simons theory \cite{AlFr}. In order to
discuss the invariants of interest in quantum gravity one has to
specify a given gauge group, $SU(2)$. It is still questionably if the
invariants that are obtained from a single gauge group are ``general
enough'' to span the whole space of states of quantum gravity. This
issue is currently under consideration.

Summarizing, we have constructed a family of ambient isotopic knot
invariants that are candidates for quantum states of the gravitational
field in the spin network representation. It is yet to be established
if this family is a complete enough ``arena'' for the quantization of
general relativity.

We wish to thank Viqar Husain and Jos\'e Labastida for pointing us to
the work of Witten and Martin. We are also grateful to Laurent Freidel
for discussions.  This work was supported in part by grants
NSF-INT-9406269, NSF-PHY-9423950, research funds of the Pennsylvania
State University, the Eberly Family research fund at PSU and PSU's
Office for Minority Faculty development. JP acknowledges support of
the Alfred P. Sloan foundation through a fellowship. We acknowledge
support of Conicyt (project 49) and PEDECIBA (Uruguay).


\begin{references}
\bibitem{Wi88} E. Witten, Commun. Math. Phys {\bf 121}, 351 (1989).
\bibitem{Kaknph} L. Kauffman ``Knots and physics'', World Scientific
Series on Knots and Everything {\bf 1}, World Scientific, Singapore (1991).
\bibitem{Gu} E. Guadagnini, ``The link invariants of the
Chern-Simons field theory, new developments in topological quantum
field theory'', De Gruyter expositions in mathematics, {\bf 10}, W. De
Gruyter, New York (1993).
\bibitem{Aw} M. Awada,  Comm. Math. Phys. {\bf 129}, 329 (1990).
\bibitem{BaNa} D. Bar-Natan, Topology, {\bf  34}, 423 (1995); q-alg/9702009.
\bibitem{BiLi} J. Birman, X. Lin, Inv. Math. {\bf 111}, 225 (1993).
\bibitem{Ba} J. Baez, Lett. Math. Phys. {\bf 26}, 43 (1992).
\bibitem{BrGaPuprl} B. Br\"ugmann, R. Gambini, J. Pullin, Phys. Rev.
Lett. {\bf 68}, 431 (1992).
\bibitem{BrGaPunpb} B. Br\"ugmann, R. Gambini, J. Pullin, Nucl. Phys.
{\bf B385}, 587 (1992).
\bibitem{BrGaPugrg} B. Br\"ugmann, R. Gambini, J. Pullin, Gen. Rel. Grav.
{\bf 25}, 1 (1993).
\bibitem{GaPubook} Gambini, R., Pullin, J.: ``Loops, knots, gauge theories and
quantum gravity.'' Cambridge: Cambridge University Press, Cambridge (1996).
\bibitem{GaTr} R. Gambini, A. Trias, Nucl. Phys. {\bf B278}, 436 (1986).
\bibitem{BrPu91} Br\"ugmann, B., Pullin, J.: Nucl. Phys. {\bf B363}, 221
(1991).
\bibitem{RoSm} C. Rovelli, L. Smolin, Nucl. Phys. {\bf B442}, 593 (1995).
\bibitem{Wi89} E. Witten, Nuc. Phys. {\bf B322}, 629 (1989).
\bibitem{GaGrPu} R. Gambini, J. Griego, J. Pullin, Phys. Lett. {\bf B} (to 
appear).
\bibitem{Ma} S. Martin, Nuc. Phys. {\bf B338}, 244 (1990).
\bibitem{KaLi} L. Kauffman, S. Lins, ``Temperley--Lieb recoupling
theory and invariants of 3-Manifolds'', Annals of Mathematics Studies,
Princeton University Press, Princeton (1994).
\bibitem{AlLa} M. Alvarez, J. M. F. Labastida, Nucl. Phys. {\bf
B433}
(1995); Erratum {\bf B441}, 403 (1995); {\bf B488}, 677 (1997);q-alg/9604010.
\bibitem{It} See for instance, C. Itzykson, J. Zuber, 
``Quantum field theory'', Mc Graw-Hill, New York (1980).
\bibitem{VaMoKh} D. Varshalovich, A. Moskalev, V. Khersonskii, 
``Quantum theory of angular momentum'', World Scientific, Singapore, (1988).
\bibitem{AsLe} A. Ashtekar, J. Lewandowski, J. Math. Phys. {\bf 5},
2170 (1995).
\bibitem{GaPucmp} R. Gambini, J. Pullin, ``Variational derivations of
exact skein relations for Chern--Simons theories'', 
Commun. Math. Phys. {\bf 185}, 621 (1997).
\bibitem{KaBa} L. Kauffman, in 
``Knots and quantum gravity'', J. Baez editor, Oxford University Press,
Oxford (1993).
\bibitem{AlFr} D. Altschuler, L. Freidel, Commun. Math. Phys. 
{\bf 170}, 41 (1995).
\end{references}
\end{document}